\documentclass[12pt,a4paper,leqno,dvips,twoside]{article}

\usepackage{graphicx}
\usepackage{amsfonts}
\usepackage{amssymb}
\usepackage{amsmath}
\usepackage{amsthm}

 \numberwithin{equation}{section}
\newtheorem{theo}{Theorem}[section]

\newtheorem{defin}{Definition}[section]

\newtheorem{cor}{Corollary}[section]
\newtheorem{prop}{Proposition}[section]
\newtheorem{rema}{Remark}[section]

\newtheorem{example}{Example}[section]

\begin{document}

\setcounter{page}{1}
\title{  Some common fixed points results on metric spaces over  topological modules}
\author{ Ion Marian Olaru }
\maketitle

\begin{abstract}
  In this paper, we replace the real numbers by a topological  R-module and define R-metric spaces $(X,d)$.
   Also, we prove some  common fixed point theorems on R-module metric
   spaces. We obtain, as a particular case the Perov theorem (see \cite{perov})
\end{abstract}

\begin{center} 2010 Mathematical Subject Classification: 47H10

Keywords: R-metric spaces, fixed point theory, topological rings, topological modules
\end{center}
\section{R-metric spaces}

In this section we shall define $R$-metric spaces and prove some properties. All axioms for an ordinary metric space can be meaningfully formulated for an abstract metric space, where the abstract metric takes  values in a partially ordered topological module of a certain type which will be defined below. Such a space will be called $R$-metric space.

We begin this section by recalling a few facts concerning topological rings, topological modules and partially ordered rings.  Unless explicitly stated otherwise all rings will be assumed to possess an identity element, denoted by $1$.
\begin{defin}(see \cite{warner}) A topology $\tau$ on a ring $(R,+,\cdot)$ is a ring topology and $R$, furnished with $\tau$, is a topological ring if the following conditions hold:

\begin{itemize}
\item[(TR 1)]$(x,y)\rightarrow x+y$ is continuous from $R\times R$ to $R$;
\item[(TR 2)]$x\rightarrow -x$ is continuous from $R$ to $R$;
\item[(TR 3)]$(x,y)\rightarrow x\cdot y$ is continuous from $R\times R$ to $R$,
\end{itemize}
where $R$ is given topology $\tau$ and $R\times R$ the cartesian product determined by topology $\tau$.
\end{defin}

\begin{defin}(see \cite{warner})
Let $R$ be a topological ring, $E$ an R-module. A topology $\mathcal{T}$ on $E$ is a R-module topology and E, furnished with $\mathcal{T}$, is a topological R-module if the following conditions hold:
\begin{itemize}
\item[(TM 1)]$(x,y)\rightarrow x+y$ is continuous from $E\times E$ to $E$;
\item[(TM 2)]$x\rightarrow -x$ is continuous from $E$ to $E$;
\item[(TM 3)]$(a,x)\rightarrow a\cdot x$ is continuous from $R\times E$ to $E$,
\end{itemize}
where $E$ is given topology $\mathcal{T}$, $E\times E$ the cartesian product determined by topology $\mathcal{T}$ and
$A\times E$ the cartesian product determined by topology of R and E.
\end{defin}

By a {\it partially ordered ring}  is meant a pair consisting of a ring and a compatible partial order, denoted by $\preceq$(see \cite{steinberg}).

In the following we always suppose that  $R$ is  an ordered topological ring such that $0\preceq 1$ and $E$  is a topological $R$-module.

\begin {defin} A subset $P$ of $E$ is called a cone if:
\begin{itemize}
\item[(i)] $P$ is closed, nonempty and $P\neq\{0_E\}$;

\item[(ii)]  $a,b\in R$, $0\preceq a$, $0\preceq b$ and $x,y\in P$ implies
$a\cdot x+b\cdot y\in P$;

\item[(iii)]$P\cap -P=\{0\}.$
\end{itemize}
\end{defin}

Given a cone $P\subset E$, we define on E the partial ordering $\leq_P$ with respect to $P$ by
\begin{equation}
x\leq_P y\  if\  and\  only\  if \ y-x\in P.
\end{equation}
We shall write $x<_Py$ to indicate that $x\leq_P y$ but $x\neq y$,
while $x\ll y$ will stand for $y-x \in int P$(interior of $P$).

\begin{example}\label{ex1}
Let $R=\mathcal{M}_{n\times n}(\mathbb{R})$ be the ring of  all matrices with n rows and n columns with entries in $\mathbb{R}$ and $E=\mathbb{R}^n$.  We define the partial order $\preceq$ on $M_{n\times n}(\mathbb{R})$ as follows

 \[A\preceq B\  if\ and\ only\ if\ for\ each\ i,j=\overline{1,n} \ we\ have\  a_{ij}\leq b_{ij}.\]

Then
\begin{itemize}
\item[(a)]the topology $\tau$, generated by matrix norm
\[N:M_{n\times n}(\mathbb{R})\rightarrow\mathbb{R},\]
\[N(A)=\max\limits_{i=\overline{1,n}}\sum\limits_{j=1}^n|a_{ij}|,\]  is a ring topology;
\item[(b)] the standard topology $\mathcal{D}$ is a R-module topology on  $\mathbb{R}^n$;
\item[(c)] $P=\{(x_1,x_2,\cdots,x_n)\in\mathbb{R}^n \mid x_i\geq 0,(\forall)i=\overline{1,n}\}$ is a cone in E.
\end{itemize}
\end{example}
Indeed, Theorem 1.3 pp 2 of \cite{warner} leads us to $(a)$.

It is obvious that $(TM\ 1)$ and $(TM\ 2)$ are satisfied. Now we consider $A_n\stackrel{\tau}{\rightarrow} A$ and $x_n\stackrel{\mathcal{D}}{\rightarrow} x$ as $n\rightarrow\infty$. Then
\[\|A_n\cdot x_n-A\cdot x\|_{\mathbb{R}^n}=\|A_n\cdot x_n-A_n\cdot x+A_n\cdot x-A\cdot x\|_{\mathbb{R}^n}\leq\]
\[\|A_n\cdot(x_n-x)\|_{\mathbb{R}^n}+\|(A_n-A)\cdot x\|_{\mathbb{R}^n}\leq N(A_n)\|x_n-x\|_{\mathbb{R}^n}+N(A_n-A)\|x\|_{\mathbb{R}^n}\stackrel{n\rightarrow\infty}{\rightarrow}0.\]
Hence, $(TM \ 3)$ holds. Thus, we have obtained $(b)$. Finally, it easy to see that $P$ is a cone in $E$.

In the following we always suppose that $E$ is a topological $R$-module, $P$ is a cone in $E$ with $intP\neq\emptyset$ and $\leq_P$ is  a partial ordering with respect to $P$.

\begin{defin}
 Let $X$ be a nonempty set. Suppose that a mapping\[d:X\times X \rightarrow  E\] satisfies:
 \begin{itemize}
 \item[$(d_1)$] $0_E\leq_P d(x,y)$ for all $x,y\in X$ and $d(x,y)=0_E$ if
 and only if $x=y$;

 \item[$(d_2)$] $d(x,y)=d(y,x)$, for all $x,y\in X$ ;

 \item[$(d_3)$] $d(x,y)\leq_P d(x,z)+d(z,y)$,  for all $x,y,z\in X$.
 \end{itemize}
  Then $d$ is called a R-metric on $X$ and $(X,d)$ is called a
  R-metric space.
 \end{defin}

\begin{example}\label{ex2}
Any cone metric space is a R-metric space.
\end{example}
\begin{example}\label{ex3}
Let $R=M_{n\times n}(\mathbb{R})$ be the ring of  all matrices with n rows and n columns with entries in $\mathbb{R}$, $E=\mathbb{R}^n$, $X=\mathbb{R}^n$ and  \[P=\{(x_1,x_2,\cdots,x_n)\in\mathbb{R}^n \mid x_i\geq 0,(\forall)i=\overline{1,n}\}\] a cone in E.

 We define the partial order $\preceq$ on $M_{n\times n}(\mathbb{R})$ as follows
 \[A\preceq B\  if\ and\ only\ if\ for\ each\ i,j=\overline{1,n} \ we\ have\ a_{ij}\leq b_{ij}.\]
 Then for all $A=(a_{ij})$, $a_{ij}>0$ we have that
 \[d:\mathbb{R}^n\times\mathbb{R}^n\rightarrow\mathbb{R}^n,\]
 \[d(x,y)=(\sum\limits_{j=1}^n a_{1j}|x_j-y_j|,\cdots,\sum\limits_{j=1}^n a_{ij}|x_j-y_j|,\cdots, \sum\limits_{j=1}^n a_{nj}|x_j-y_j|)\]
 is a $R$-metric on $X$.
\end{example}
Indeed,
\begin{itemize}
\item[$(d_1)$] Since $\sum\limits_{j=1}^n a_{ij}|x_j-y_j|\geq 0$ for all $i=\overline{1,n}$, we have that $0\leq_P d(x,y)$ for all $x,y\in \mathbb{R}^n$. Also, $d(x,y)=0$ involve that $\sum\limits_{j=1}^n a_{ij}|x_j-y_j|=0$ which means that $x_j=y_j$ for all $j=\overline{1,n}$. It follows that $x=y$.
\item[$(d_2)$] It is obvious that $d(x,y)=d(y,x)$, for all $x,y\in \mathbb{R}^n$.
\item[$(d_3)$]  Let be $x,y,z\in \mathbb{R}^n$. Since $\sum\limits_{j=1}^n a_{ij}|x_j-y_j|\leq \sum\limits_{j=1}^n a_{ij}|x_j-z_j|+\sum\limits_{j=1}^n a_{ij}|z_j-y_j|$, we have  that $d(x,y)\leq_P d(x,z)+d(z,y)$,  for all $x,y,z\in \mathbb{R}^n$.
\end{itemize}
\pagebreak
In the following, we shall write $x \prec y$ to indicate that $x \preceq y$ but $x\neq y$.
\begin{rema}\label{r1} We have that:
\begin{itemize}
\item[(i)]$int P+ int P\subseteq int P$;
\item[(ii)]$\lambda\cdot int P\subseteq int P$, where $\lambda$ is a invertible element of the ring $R$ such that $0\prec\lambda$;
\item[(iii)] if  $x\leq_P y $ and $0 \preceq \alpha$, then $\alpha\cdot x\preceq\alpha\cdot y$.
\end{itemize}\end{rema}

\noindent{\bf Proof:}

$(i)$ Let  be $x\in int P+int P$. Then there exists $x_1,x_2 \in int P$ such that $x=x_1+x_2$. It follows that there exists $V_1$ neighborhood of $x_1$ and  $V_2$ neighborhood of $x_2$ such that

\[x_1\in V_1\subset P \ and\  x_2\in V_2\subset P.\]

Since for each $x_0\in E$, the mapping $x\rightarrow x+x_0$ is a
homeomorphism of $E$ onto itself,we have that $V_1+V_2$ is a neighborhood of $x$ with respect to topology $\mathcal{T}$. Thus, $x\in int P$.

$(ii)$ Let $0\prec\lambda$ be an invertible element of  the ring $R$ and $x=\lambda\cdot c$, $c\in int P$. It follows that there exists a neighborhood $V$ of $c$ such that $c\in V\subset P$.

Since the mapping $x\rightarrow \lambda x$ is a homeomorphism of $E$ onto itself, we have that $\lambda\cdot V$ is a neighborhood of $x$ with respect to the topology $\mathcal{T}$. Thus, $x\in int P$.

$(iii)$  If  $x\leq_P y $,  then $y-x\in P$. It follows that for all $0 \preceq \alpha$ we have that $\alpha\cdot (y-x)\in P$ i.e. $\alpha\cdot x\preceq\alpha\cdot y$.

In the following, unless otherwise specified, we always suppose
that there exists  at least one  a sequence $\{\alpha_n\}\subset
R$ of invertible elements such that  $0\prec\alpha_n$  and
$\alpha_n\rightarrow 0$ as $n\rightarrow \infty$.

\begin{rema}\label{r2} Let $E$ be a $R$-topological module and $P\subset E$ a cone. We have that:
\begin{itemize}
\item[(i)] If $u\leq_P v$ and $v\ll w$, then $u\ll w$;
\item[(ii)]If $u\ll v$ and $v\leq_P w$, then $u\ll w$;
\item[(iii)]If $u\ll v$ and $v\ll w$, then $u\ll w$;

\item[(iv)] If $0\leq_P u\ll c$ for each $c\in Int P$, then $u=0$;
\item[(v)] If $a\leq_P b+c $ for each $c\in Int P$, then $a\leq_P
b$;

\item[(vi)]If $0\ll c$, $0\leq_P a_n$ and $a_n\rightarrow 0$ then
there exists $n_0\in \mathbb{N}$ such that $a_n\ll c$ for all
$n\geq n_0$.

\end{itemize}
\end{rema}
\noindent{\bf Proof:}
\begin{itemize}
\item[(i)]   We have to prove that $w - u \in int P$ if $v - u \in
P$ and $w - v\in int P$.   The condition $(TM_1)$ implies that there
exists a neighborhood $V$ of $0$ such that $w-v+V\subset P$. It
follows that $w-u+V=(w-v)+V+(v-u)\subset P+P\subset P$. Since for
each $x_0\in E$ the mapping $x\rightarrow x+x_0$ is a
homeomorphism of $E$ onto itself we have that $w-u+V$ is a neighborhood of $w-u$ with respect to the topology $\mathcal{T}$. Thus, $w-u\in int\  P$.
\item[(ii)] Analogous with $(i)$.
\item[(iii)] We have to prove that $w - u \in int P$ if $v - u \in
int P$ and $w - v\in int P$.  Since we have $int P+int P\subset int P$, it easy to see that the above assertion  is satisfied.

\item[(iv)]  Let $\{\alpha_n\}_{n\in\mathbb{N}}\subset R$ be a sequence of invertible elements
such that   for each $n\in \mathbb{N}$  we have $0\prec\alpha_n$  and $\alpha_n\rightarrow 0$ as $n\rightarrow \infty$.

Since for each invertible element $\lambda_0\in R$  the mapping
$x\rightarrow \lambda_0\cdot x$ is a homeomorphism of $E$ onto
itself, we have  that  if  $V$  is a neighborhood of zero then  $ \lambda_0\cdot V$ is a neighborhood of zero. Hence,
$\alpha_n\cdot c\in int P$ for each $n\in \mathbb{N}$.

Then $\alpha_n\cdot c-u\in int P$. It follows that
$\lim\limits_{n\rightarrow\infty}\alpha_n\cdot
c-u=-u\in\overline{P}=P$. Thus, $u\in P\cap -P=\{0\}$. \item[(v)]
Analogous with $(iv)$.

\item[(vi)] Let be $0\ll c$, $0\leq_P a_n$ and $a_n\rightarrow 0$.
Since for each invertible element $\lambda_0\in R$ the mapping
$x\rightarrow \lambda_0\cdot x$ is a homeomorphism of $E$ onto
itself, it follows that  for all neighborhood $V$ of zero we have
that $- V$ is a neighborhood of zero. Now, $0\ll c$ implies that
there exists a neighborhood  $U $ of zero such that $c+U\subset
P$. Let $V=U\cap -U$  be a neighborhood of zero. Since $a_n $
converges to zero, it follows that there exists $n_0\in\mathbb{N}$
such that $a_n\in V$ for all $n\geq n_0$. Then  we have that
$c-a_n\in c+V\subset c+U\subset P$ for all $n\geq n_0$. Thus,
$a_n\ll c$ for all $n\geq n_0$.

\end{itemize}
\begin{defin}\label{def2}
 Let $\{x_n\}$ be a sequence in a  R-metric space $(X,d)$ and $x\in X$. We say that:
\begin{itemize}

\item[(i)] the sequence $\{x_n\}$ converges to $x$ and is denoted by
$\lim\limits_{n\rightarrow\infty}x_n=x$ if for every $0\ll c$ there exists
$N\in\mathbb{N}$ such that $d(x_n,x)\ll c$, for all $n>N$;
\item[(ii)]  the sequence $\{x_n\}$ is  a Cauchy sequence if for every $c\in E$, $0\ll c$ there exists $N\in\mathbb{N}$ such that $d(x_m,x_n)\ll c$ for all $m,n>N$;
\end{itemize}
The  R-metric space $(X,d)$ is called complete if every Cauchy sequence is convergent.
\end{defin}

From the above remark we obtain
\begin{rema}\label{r3}Let $(X,d)$ be a  R-metric space and $\{x_n\}$ be a sequence in $X$. If $\{x_n\}$ converges to $x$ and $\{x_n\}$ converges to $y$, then $x = y$.\end{rema}
Indeed, for all $0\ll c$
\[d(x,y)\leq_P d(x,x_n)+d(x_n,y)\ll 2\cdot c.\]
 Hence, $d(x,y)=0$ i.e. $x=y$.
\section{Common fixed points theorems}
In  this section we obtain several coincidence and common fixed point  theorems for mappings defined on a $R$-metric space.
\begin{defin}(see \cite{abas}) Let $f$ and $g$ be self maps of a set $X$. If $w = fx = gx$ for some
$x\in X$, then $x$ is called a coincidence point of $f$ and $g$,
and $w$ is called a point of coincidence of $f$ and
$g$.\end{defin}
Jungck \cite{jungck}, defined a pair of self mappings to be weakly compatible if they commute at their coincidence points.
\begin{prop}\label{p1}(see \cite{abas}) Let $f$ and $g$ be weakly compatible self maps of a set $X$. If $f$ and $g$ have a unique point of coincidence
$w = fx = gx$, then $w$ is the unique common fixed point of $f$
and $g$.
\end{prop}
Let $\mathcal{K}$ be the set of all $k\in R$, $0\preceq k$ which have the property that there exists  a unique $S\in R$ such that $S =\lim\limits_{n\rightarrow\infty}(1+k+\cdots+k^n)$.
\begin{example}\label{ex4}
Let be  $A\in\mathcal{M}_{n\times n}(\mathbb{R}_+)$ such that $\rho(A)<1$. Then $A\in \mathcal{K}$.
\end{example}

\begin{theo}\label{t1}
Let $(X,d)$ be a R-metric space and suppose  that the mappings $f,g:X\rightarrow X$ satisfy:
\begin{itemize}
\item[(i)]the range of $g$ contains the range of $f$ and $g(X)$ is a complete subspace of $X$;
\item[(ii)] there exists $k\in \mathcal{K}$ such that $d(fx,fy)\leq_P  k\cdot d(gx,gy)$ for all $x,y\in X$.
\end{itemize}
Then $f$ and $g$ have a unique point of coincidence in $X$. Moreover, if $f$ and $g$ are weakly compatible then, $f$ and $g$ have a unique common fixed point.
\end{theo}
\noindent{\bf Proof:} Let $x_0$ be an arbitrary point in $X$. We
choose a point $x_1\in X$ such that $f(x_0)=g(x_1).$ Continuing
this process, having chosen $x_n\in X$, we obtain $x_{n+1}\in X$
such that $f(x_n)=g(x_{n+1})$. Then
\[d(gx_{n+1},gx_n)=d(fx_n,fx_{n-1})\leq _P k \cdot d(gx_n, gx_{n-1})\leq_P \]
\[\leq_P k^2 \cdot d(gx_{n-1}, gx_{n-2})\leq_P\cdots\leq_P k^n \cdot d(gx_1,gx_0).\]
 We  denote $S_n=1+k+\cdots+k^n$and we get that
 \[d(gx_n,gx_{n+p})\leq_P  d(gx_n,gx_{n+1})+d(gx_{n+1},gx_{n+2})+\cdots+d(gx_{n+p-1},gx_{n+p})\leq_P\]
 \[\leq_P (k^n+ k^{n+1}+\cdots+ k^{n+p-1})\cdot d(gx_1,gx_0)=(S_{n+p-1}-S_{n-1})\cdot d(gx_1,gx_0) ,\] for all $p\geq 1$.
 Thus, via Remark \ref{r2} $(vi)$, we obtain that  $gx_n$ is a Cauchy
 sequence. Since $g(X)$ is complete, there exists $q\in g(X)$ such that $gx_n\rightarrow q $ as $n\rightarrow \infty$.
 Consequently, we can find $p\in X$ such that $g(p)=q$. Further, for each $0\ll c$ there exists $n_0\in\mathbb{N}$
 such that for all $n\geq n_0$
 \[d(gx_n,fp)=d(fx_{n-1},f(p))\leq_P k \cdot d(gx_{n-1},gp)\ll c.\]
 It follows that $gx_n\rightarrow fp$ as $n\rightarrow\infty$. The uniqueness of the limit implies that $fp=gp=q$.
 Next we show that $f$ and $g$ have a unique point of coincidence.
 For this, we assume that there exists another point $p_1\in X$ such that $fp_1=gp_1$.
 We have
 \[d(gp_1,gp)=d(fp_1,fp)\leq_P k\cdot d(gp_1,gp)=k\cdot d(fp_1,fp)\leq_P k^2\cdot d(gp_1,gp)\leq_P\cdots\leq_P k^n\cdot d(gp_1,gp).\]
 Let  be $0\ll c$. Since $k^n\rightarrow 0$ as $n\rightarrow\infty$ it follows that there
 exists $n_0\in \mathbb{N}$ such that $k^n\cdot d(gp_1,gp)\ll c$ for all $n\geq n_0$.
 Then $d(gp_1,gp)\ll c$ for each $0\ll c$. Thus $d(gp_1,gp)=0$
 i.e. $gp_1=gp$. From Proposition \ref{p1} $f$ and $g$ have a unique common fixed point.
 \begin{rema}\label{r5} The above theorem generalizes Theorem 2.1 of Abbas and Jungck \cite{abas}, which itself is a generalization of  Banach fixed point theorem.\end{rema}

 \begin{cor}\label{c1} Let $(X,d)$ be a  complete  $R$-metric space and we suppose that the mapping $f:X\rightarrow X$ satisfies:
\begin{itemize}
\item[(i)] there exists $k\in \mathcal{K}$ such that $d(fx,fy)\leq_P  k\cdot
d(x,y)$ for all $x,y\in X$.
\end{itemize}
Then $f$ has in $X$ a unique  fixed point point.
\end{cor}
\noindent{\bf Proof:} The proof uses Theorem \ref{t1} by replacing
$g$ with   the identity mapping.

From the above corollary  using  Example \ref{ex1}, we obtain the
Perov fixed point theorem (see \cite{perov})
\begin{cor}
Let $(X,d)$ be a  complete $\mathcal{M}_{n\times n}(\mathbb{R}_+)-$ metric space and $E=\mathbb{R}^n$ and we suppose  the mapping $f:X\rightarrow X$ satisfies:
\begin{itemize}
\item[(i)] there exists $A\in\mathcal{M}_{n\times n}(\mathbb{R}_+)$ with  $\rho(A)<1$ such that

\[d(fx,fy)\leq_P  A\cdot d(x,y),\]for all $x,y\in X$.
\end{itemize}
Then $f$ has in $X$ a unique  fixed point point.
\end{cor}
\section{Comparison function}
\begin{defin}\label{def3} Let P be a cone in a topological R-module E. A function $\varphi:P \rightarrow P$ is called a comparison function if
\begin{itemize}
\item[(i)]$\varphi(0)=0$ and $ \varphi(t)<_P t$ for all  $t\in
P-\{0\}$;

\item[(ii)]$t_1\leq_P t_2$ implies
$\varphi(t_1)\leq_P\varphi(t_2)$; \item[(iii)] $t\in int P$
implies $t-\varphi(t)\in  int P$; \item[(iv)] if $t\in P-\{0\}$
and $0\ll c$, then there exists $n_0\in\mathbb{N}$ such that
$\varphi^n(t)\ll c$ for each $n\geq n_0$.
\end{itemize}\end{defin}
\begin{example}\label{ex5}
Let  P be a cone in a topological R-module E and $\lambda\in \mathcal{K}$ such that $0\prec 1-\lambda$. Then  $\varphi:P\rightarrow P$, defined by $\varphi(t)=\lambda \cdot t$ is a comparison function.
\end{example}
Indeed,

$(i)$ It is obvious that $\varphi(0)=0$ and $ \varphi(t)<_P t$ for
all  $t\in P-\{0\}$.

$(ii)$ if $t_1\leq_P t_2$  and $\lambda\in \mathcal{K}$ then $\lambda\cdot(t_2-t_1)\in P$. Thus $\varphi(t_1)\leq_P\varphi(t_2)$.

$(iii)$  We remark that if $\lambda\in \mathcal{K}$, then $1-\lambda$ is an invertible element of  the ring $R$. Now, let be $t\in int P$.  Then $(1-\lambda)\cdot t\in (1-\lambda) int P\subset int P$.

$(iv)$ Let be  $t\in P-\{0\}$ and $0\ll c$. Then
\[\varphi^n(t)=\lambda^n\cdot t\stackrel{n\rightarrow\infty} {\rightarrow}0.\]
We obtain, via Remark \ref{r2}, that there exists $n_0\in\mathbb{N}$ such that $\varphi^n(t)\ll c$ for each $n\geq n_0$.

Let $(X,d)$ be a $R$-metric space and let $\varphi:K\rightarrow K$ be a comparison function. For a pair $(f,g)$ of self-mappings on $X$ consider the following condition
\begin{itemize}
\item[(C)] for arbitrary $x,y\in X$ there exists $u\in \{d(gx,gy),d(gx,fx), d(gy,fy)\}$ such that $d(fx,fy)\leq_P \varphi(u)$.
\end{itemize}
\begin{theo}\label{t1} Let $(X,d)$ be  a R-metric space and let $f,g:X\rightarrow X$  such that
\begin{itemize}
\item[(i)]   the pair $(f,g)$ satisfies the  condition (C) for some comparison function $\varphi$;
\item[(ii)] $f(X)\subset g(X)$ and f(X) or g(X) is a complete subspace of X.
\end{itemize} Then f and g have a unique point of coincidence in X. Moreover if $f$ and $g$ are weakly compatible, then $f$ and $g$ have a unique common fixed point. \end{theo}
\noindent{\bf Proof:} Let $x_0$ be an arbitrary point in $X$. We choose a point $x_1\in X$ such that $fx_0=gx_1.$ Continuing this process, having chosen $x_n\in X$, we obtain $x_{n+1}\in X$ such that $fx_n=gx_{n+1}$.

 $\ulcorner$ We shall prove that the sequence $\{y_n\}$, where $y_n=fx_n=gx_{n+1}$( the so-called Jungck sequence ) is a Cauchy sequence in $R$-metric space $(X,d)$.

 If $y_N=y_{N+1}$ for some $N\in\mathbb{N}$ then $y_m=y_N$ for each $m>N$ and the conclusion follows. Indeed, we prove by induction  arguments that

 \begin{equation}\label{e4}y_{N+k}=y_{N+k+1}, (\forall) k\in\mathbb{N}.\end{equation}
 For $k=0$ we have $y_N=y_{N+1}$. Let (\ref{e4}) hold for all $k=\overline{0,i}$. Then
 \[d(y_{N+i+1},y_{N+i+2})= d(fx_{N+i},fx_{N+i+1})\leq_P\varphi(u),\]
 where
 \[u\in\{d(gx_{N+i},gx_{N+i+1}), d(gx_{N+i},fx_{N+i}), d(gx_{N+i+1},fx_{N+i+1})\}=\]
 \[\{d(y_{N+i-1},y_{N+i}),d(y_{N+i-1},y_{N+i}),d(y_{N+i},y_{N+i+1})\}=\{0\}.\]
 Hence, $d(y_{N+i+1},y_{N+i+2})\leq_P\varphi(u)=0 $ i.e. $y_{N+i+1}=y_{N+i+2}$. Q.E.D

 Suppose that $y_n\neq y_{n+1}$ for each $n\in\mathbb{N}$.  The condition (C) implies that
\[d(y_n,y_{n+1})=d(fx_n,fx_{n+1})\leq_P\varphi(u),\]  where
\[u\in \{d(gx_n,gx_{n+1}), d(fx_n,gx_n), d(fx_{n+1},gx_{n+1})\}=\{d(y_{n-1},y_n), d(y_n,y_{n-1}), d(y_{n+1},y_n)\}.\]
The case $u=d(y_{n+1},y_n)$ is impossible, since this would imply
\[d(y_{n+1},y_n)\leq_P \varphi(d(y_{n+1},y_n))<_P d(y_{n+1},y_n).\]
Thus, $u=d(y_n,y_{n-1})$ and
 \[d(y_{n+1}, y_n)\leq_P \varphi(d(y_n,y_{n-1}))\leq_P\cdots\leq_P \varphi^n(d(y_1,y_0)).\]
Using Remark \ref{r2} (i) and property $(iv)$ of  the comparison function we obtain that  for all $0\ll\varepsilon$ there exists $n_0\in \mathbb{N}$ such that
\begin{equation}\label{e1}d(y_n,y_{n+1})\ll \varepsilon,\ (\forall) n\geq n_0.\end{equation}
Now, let  be $0\ll c$.
 Then, using property $(iii)$ of the comparison function, we get that
\begin{equation}\label{e2}d(y_n,y_{n+1})\ll c-\varphi(c),\ (\forall) n\geq n_0.\end{equation}
Let us fix now $n\geq n_0$ and let us prove that
\begin{equation}\label{e3}d(y_n,y_{k+1})\ll  c,\ (\forall) k\geq n.\end{equation}
Indeed, for $k=n$ we have
\[d(y_n,y_{n+1})\ll c-\varphi(c)\leq_P c.\]
Hence,
\[d(y_n,y_{n+1})\ll c.\]
Let (\ref{e3} ) hold for some $k\geq n$. Then we have
\[d(y_n,y_{k+2})\leq_P d(y_n,y_{n+1})+ d(y_{n+1},y_{k+2})\ll \]
\[c-\varphi(c)+d(fx_{n+1},fx_{k+2})\leq_P c-\varphi(c)+\varphi(u),\] where
\[u\in\{d(gx_{n+1},gx_{k+2}),d(gx_{n+1},fx_{n+1}),d(gx_{k+2},fx_{k+2})\}.\]
Consider now the following three possible cases:
\begin{itemize}
\item[ \it Case 1:]$ u=d(gx_{n+1},gx_{k+2})$. Then
\[\varphi(u)=\varphi(d(gx_{n+1},gx_{k+2}))=\varphi(d(y_n,y_{k+1}))\leq_P \varphi(c).\]
From the above relation it follows that,
\[d(y_n,y_{k+2})\ll c-\varphi(c)+\varphi(u)\leq_P c-\varphi(c)+\varphi(c)=c.\]
Hence, $d(y_n,y_{k+2})\ll c.$
\item[\it Case 2:] $u=d(gx_{n+1},fx_{n+1})=d(y_n,y_{n+1})$. Then
\[\varphi(u)\leq_P \varphi(d(y_n,y_{n+1}))\leq_P \varphi(c-\varphi(c))\leq_P \varphi(c).\]
From the above relation we get that,
\[d(y_n,y_{k+2})\ll c-\varphi(c)+\varphi(u)\leq_P c-\varphi(c)+\varphi(c)=c.\]
Hence, $d(y_n,y_{k+2})\ll c.$
\item[\it Case 3:] $u=d(gx_{k+2},fx_{k+2})$. Then
\[\varphi(u)=\varphi(d(gx_{k+2},fx_{k+2}))=\varphi(d(y_{k+1},y_{k+2}))\leq_P\varphi(c-\varphi(c))\leq_P \varphi(c).\]
From the above relation we get that,
\[d(y_n,y_{k+2})\ll c-\varphi(c)+\varphi(u)\leq_P c-\varphi(c)+\varphi(c)=c.\]
Hence, $d(y_n,y_{k+2})\ll c.$
\end{itemize}
So, it has been proved by induction that  $\{y_n\}$ is a Cauchy sequence.$\lrcorner$

Since, by assumption, $f(X)$ or $g(X)$ is a complete subspace of
$X$, we conclude that there exists $q\in g(X)$ such that
$y_n=fx_n=gx_{n+1}\rightarrow q$ as $n\rightarrow\infty$ and there
exists $p\in X$ such that $q=gp$.

$\ulcorner$ We claim that $q=fp$. Indeed, if we suppose that $d(fp,q)\neq 0$,  then we have
\[d(fp,q)\leq_P d(fp,fx_n)+ d(fx_n,q)\leq_P \varphi(u)+ d(y_n,q),\]
where
\[u\in\{d(gp,gx_n), d(gp,fp),d(gx_n,fx_n)\}\]
Let $0\ll c$. At least one of the following three cases holds for infinitely many $n\in \mathbb{N}$:
\begin{itemize}
\item[\it Case 1:]$u=d(gp,gx_n)$. Then, there exists $n_0(c)\in\mathbb{N}$ such that for all $n\geq n_0(c)$
\[d(fp,q)\leq_P\varphi(d(gp,gx_n))+ d(y_n,q)<_P  d(q,y_{n-1})+ d(y_n,q)\ll 2\cdot c.\]
It follows that $d(fp,q)=0$, which is a contradiction.
\item[\it Case 2:] $u=d(gp,fp)=d(fp,q)$ . Then we have
\[d(fp,q)\leq_P \varphi (d(q,fp))+ d(y_n,q)\ll  \varphi (d(q,fp))+c.\]
Thus,  $d(fp,q)\leq_P \varphi (d(q,fp))$. So, by using  of the properties $(ii)$ and $(iv)$ of the comparison function, we obtain that there exists $n_0\in\mathbb{N}$ such that $d(fp,q)\leq_P \varphi^n(d(fp,q))\ll c $ i.e. $d(fp,q)=0$, which is a contradiction.
\item[\it Case 3:] $u=d(gx_n,fx_n)=d(y_{n-1},y_n)$. Then, there exists $n_0(c)\in\mathbb{N}$ such that for all $n\geq n_0(c)$ we have
     \[d(fp,q)\leq_P\varphi (d(y_{n-1},y_n))+ d(y_n,q)<_P d(y_{n-1},y_n)+ d(y_n,q)\ll 2\cdot c,\] i.e.  $d(fp,q)=0$, which is a contradiction.
\end{itemize}
It follows that $fp=gp=q$ i.e.  $p$ is a coincidence point of the pair $(f,g)$ and $q$ is a point of coincidence.$\lrcorner$

$\ulcorner$ Next we show that $f$ and $g$ have a unique point of coincidence. For this we assume that there exists another point $p_1\in X$ such that $fp_1=gp_1$. If we suppose that $d(fp_1,fp)\neq 0$  we get that $d(fp_1,fp)\leq_P \varphi(u)$, where
\[u\in \{d(gp_1,gp), d(gp_1,fp_1),d(gp,fp)\}=\{d(gp_1,gp), 0\}.\]
In both possible cases a contradiction  easily follows :
$d(fp_1,fp)\leq_P \varphi(d(fp_1,fp))<_P d(fp_1,fp)$ or
$d(fp_1,fp)\leq_P \varphi(0)=0$. We conclude that the mappings $f$
and $g$ have a unique point of coincidence. From Proposition
\ref{p1} $f$ and $g$ have a unique common fixed point.$\lrcorner$

\noindent Department of Mathematics,\\
Faculty of Science,\\
University "Lucian Blaga" of Sibiu,\\
Dr. Ion Ratiu 5-7, Sibiu, 550012, Romania \\
E-mail: marian.olaru@ulbsibiu.ro

\end{document}